\documentclass[12pt]{amsart}
\usepackage{amsmath,amsthm,amsfonts,amssymb}
\begin{document} 
\newcommand{\B}{{\mathbb B}}
\newcommand{\C}{{\mathbb C}}
\newcommand{\N}{{\mathbb N}}
\newcommand{\Q}{{\mathbb Q}}
\newcommand{\Zeta}{\mathop{\Z'}}
\newcommand{\Z}{{\mathbb Z}}
\renewcommand{\P}{{\mathbb P}}
\newcommand{\R}{{\mathbb R}}
\newcommand{\rc}{\subset}
\newcommand{\rank}{\mathop{rank}}
\newcommand{\trace}{\mathop{tr}}
\newcommand{\dimc}{\mathop{dim}_{\C}}
\newcommand{\Lie}{\mathop{Lie}}
\newcommand{\Auto}{\mathop{{\rm Aut}_{\mathcal O}}}
\newcommand{\alg}[1]{{\mathbf #1}}
\newcommand{\tensor}{\otimes}
\newtheorem*{definition}{Definition}
\newtheorem*{claim}{Claim}
\newtheorem{corollary}{Corollary}
\newtheorem*{Conjecture}{Conjecture}
\newtheorem*{SpecAss}{Special Assumptions}
\newtheorem*{example}{Example}
\newtheorem*{remark}{Remark}
\newtheorem*{observation}{Observation}
\newtheorem*{fact}{Fact}
\newtheorem*{remarks}{Remarks}
\newtheorem{lemma}{Lemma}
\newtheorem{proposition}{Proposition}
\newtheorem*{theorem}{Theorem}
\title{%
A splitting criterion for two-dimensional semi-tori
}
\author {J\"org Winkelmann}
\begin{abstract}
We investigate conditions under which a two-dimensional
complex semi-torus splits into a direct product of $\C^*$
and a one-dimensional compact complex torus.
\end{abstract}
%
\address{%
J\"org Winkelmann \\
 Institut Elie Cartan (Math\'ematiques)\\
 Universit\'e Henri Poincar\'e Nancy 1\\
 B.P. 239\\
 F-54506 Vand\oe uvre-les-Nancy Cedex\\
 France
}
\email{jwinkel@member.ams.org\newline\indent{\itshape Webpage: }%
http://www.math.unibas.ch/\~{ }winkel/
}
\thanks{
{\em Acknowledgement.}
The author wants to thank 
the Korea Institute for Advanced Study in Seoul.
The research for this article was done during the stays of the
author at these institutes.}
\maketitle
A semi-torus is a complex Lie group arising as a quotient of the
additive group of a complex vector space $V$ by a discrete subgroup
$\Gamma$ with the property that $\Gamma$ generates $V$ as complex
vector space.
Semi-tori without non-constant holomorphic functions are also
known as {\em Cousin groups}. They have been studied by many
mathematicians since the work of Cousin about a century ago
\cite{C}.

In this article we investigate a special aspect of two-dimensional
semi-tori.
Let $\Gamma$ be a discrete subgroup of $\Z$-rank three in $(\C^2,+)$.
Then $T=\C^2/\Gamma$ is a {\em semi-torus}. Let $\Gamma_\R$ denote
the $\R$-span of $\Gamma$ and $H=\Gamma_\R\cap i\Gamma_\R$.
Then for every $\gamma\in\Gamma\setminus H$ the quotient 
$S_\gamma=\C\gamma/(\C\gamma\cap\Gamma)$ is a closed complex 
Lie subgroup of $T$ which is isomorphic to
$\C^*$. 
The quotient $E_\gamma=T/S_\gamma$ is an elliptic curve.
Usually, for different elements $\gamma\in\Gamma$ the isomorphism class
and even the isogeny class of this elliptic 
curve $E_\gamma$ will depend on $\gamma$.
For example, if $T$ is the quotient of $\C^2$ by the lattice 
\[
\Gamma=\left< 
\begin{pmatrix} 1  \\ 0 \end{pmatrix},
\begin{pmatrix} 0 \\ 1 \end{pmatrix},
\begin{pmatrix} \tau \\ \sigma \end{pmatrix}
\right>_\Z,
\]
with $\tau,\sigma\in H^+=\{z\in\C:\Im(z)>0\}$,
then there is a  quotient elliptic curve isomorphic to
$\C/\left<1,\sigma\right>$ as well as one isomorphic to
$\C/\left<1,\tau\right>$.
On the other hand, if $T\simeq\C^*\times E$ for some elliptic curve
$E$, then the restriction of the projection map $T\to E_\gamma$
to $\{1\}\times E$ yields an isogeny between $E$ and $E_\gamma$.
Thus in this case all these quotient elliptic curves $E_\gamma$ must
be isogenous. 
In this paper we investigate inhowfar this property characterizes
those $T$ which split into a direct product of $\C^*$ and an elliptic curve.
It turns out that besides such a splitting also a certain arithmetic
property may cause all the quotient elliptic curves to be isogenous.

\begin{theorem}
Let $\Gamma$ be a discrete subgroup of $\Z$-rank three in $(\C^2,+)$
and $T=\C^2/\Gamma$.
Assume that for every surjective complex Lie group homomorphism
from $T$ to an elliptic curve $E$ this elliptic curve is of the same
isogeny class.

Then one of the following conditions hold:
\begin{enumerate}
\item
$T$ is isomorphic to a direct product of ${}\,\C^*$ and an elliptic curve, or
\item
There is a number field $k$ of degree $3$ over $\Q$ such that
$\Gamma\subset k^2$ after some linear change of coordinates on $\C^2$.
\end{enumerate}
Conversely, if one of these two conditions is fulfilled, then there
exists an elliptic curve $E_0$ such that $E$ is isogenous to $E_0$ 
for every surjective complex Lie group homomorphism from $T$ to an
elliptic curve $E$.
\end{theorem}

It should be remarked that $T=\C^2/\Gamma$ must be a direct product
of $\C^*$ and an elliptic curve
if $\Gamma\subset k^2$ for some {\em quadratic} number field
(see lemma~\ref{q-field}).

One might inquire what happens if one asks for isomorphisms instead
of mere isogenies between all quotient elliptic curves.
As it turns out this simply is too much to ask for:
There are always non-isomorphic quotient elliptic curves,
as we will see in the last section.
\section{Proofs}
Before we prove the theorem, we need to deduce some auxiliary results.
\begin{lemma}\label{q-field}
Let $k$ be a quadratic number field, ${\mathcal O}_k$ its ring of algebraic
integers and $\Gamma\subset k^2$ a
subgroup of $\Z$-rank three which is discrete in $\C^2$.

Then there exists an elliptic curve $E$ (isogenous to $\C/{\mathcal O}_k$)
such that $\C^2/\Gamma\simeq\C^*\times E$.
\end{lemma}
\begin{proof}
Since $\Gamma$ is discrete, $k$ can not be totally real. Hence
$k=\Q[\sqrt{-n}]$ for some $n\in\N$. Let $\tau=\sqrt{-n}\in k$.
Let $\Gamma_\Q=\Gamma\tensor_\Z\Q$. Consider 
$H=\Gamma_\Q\cap \tau(\Gamma_\Q)$.
Since $H$ is the intersection of two $\Q$-hyperplanes in $k^2$, it is
clear that $\dim_\Q(H)=2$. 
Moreover $H$ is stable under scalar multiplication with elements of $k$.
Hence $H\tensor_\Q\R$ is a complex line which we will
call $H_{\C}$.
The intersection $H\cap\Gamma$ is a lattice of
$\Z$-rank $2$ in $H_{\C}$. Since $H$ is $k$-vector space, and since
${\mathcal O}_k$ is a lattice of $\Z$-rank 2, it follows that
$H_\C/(H\cap\Gamma)$ is isogenous to $\C/{\mathcal O}_k$.
On the other hand, $\rank_\Z(H\cap\Gamma)=2$ implies that
$\rank_{\Z}\pi(\Gamma)=1$ where $\pi$ denotes the natural
projection $\C^2\to \C^2/H_\C$. Therefore $\pi(\Gamma)\simeq\Z$ and
$\Gamma$ contains an element $\gamma_0$
such that $\pi(\gamma_0)$ generates $\pi(\Gamma)$.
It follows that 
\[
\C^2/\Gamma\simeq (\C\gamma_0/\Z\gamma_0)\times 
(H_\C/(H\cap\Gamma)).
\] 
This implies the statement.
\end{proof}

\begin{lemma}
Let $k$ be a cubic number field. Then $GL_2(\Q)$ has exactly
two orbits in $\P_1(k)$, namely $\P_1(\Q)$ and its complement.
\end{lemma}
\begin{proof}
We start by the claim: {\em $PGL_2(\Q)$ acts freely on 
$\Omega=\P_1(k)\setminus\P_1(\Q)$.}

Indeed, let $x\in k\setminus\Q$. Then $\left(\begin{smallmatrix} a & b \\
c & d 				       \end{smallmatrix}\right)
$ is in the isotropy at $[x:1]$ iff
\[
x(cx+d)=ax+b \iff cx^2 + (d-a)x -b = 0.
\]
Now $x\in k\setminus\Q$ and $\deg(k/\Q)=3$ imply that $1,x,x^2$
are all $\Q$-linearly independent.
Thus the above equation implies that $c=b=0$ and $a=d$.
This yields the claim.

Next we consider the Borel group 
\[
B=\left\{ \begin{pmatrix} a & b \\ & d \end{pmatrix} : a,d\in\Q^*, b\in\Q
\right\}.
\]
Such an element in $B$ maps $[x:1]$
to $[ax+b:d]=[\frac{a}{d}x+\frac{b}{d}:1]$.
Fix a generator (``primitive element'') 
$\tau$ for the field extension  $k/\Q$.
Then every element of $k$ be written uniquely as $a+b\tau+c\tau^2$
with $a,b,c\in\Q$.
We define a map $\phi:\Omega\to\P_1(\Q)$ as follows:
\[
[a+b\tau+c\tau^2:1] \to [b:c].
\]
The fibers of $\phi$ are precisely the $B$-orbits in $\Omega$.

We continue with a second claim: {\em There exists a cubic polynomial $Q\in\Q[X]$
such that for all $\lambda\in\Q$ with $Q(\lambda)\ne 0$
there exists an element $A\in GL_2(\Q)$ such that
$[\lambda\tau+\tau^2:1]=A([\tau:1])$}.

Let $\tau^3=p_2\tau^2+p_1\tau+p_0$ ($p_i\in\Q)$).
Choose $c=1$, $d=-(\lambda+p_2)$, $a=p_1-\lambda(\lambda+p_2)$
and $b=p_0$.

Then
\[
\det\begin{pmatrix} a & b \\ c & d \end{pmatrix}
= -p_1(\lambda+p_2) + \lambda(\lambda+p_2)^2-p_0
\]
and
\[
(\lambda\tau+\tau^2)(c\tau+d)=(a\tau+b).
\]
This proves the second claim
with
$A=\left( \begin{smallmatrix} a & b \\ c & d \end{smallmatrix} \right)$
and
$Q(\lambda)=  -p_1(\lambda+p_2) + \lambda(\lambda+p_2)^2-p_0$.

Now we can prove the lemma. Let $G=PGL_2(\Q)$.
By the last claim, the $G$-orbit through $[\tau:1]$ contains
all of $\Omega$ with the possible exception of at most three
$B$-orbits. But $G$ acts freely on $\Omega$ and $G/B\simeq\P_1(\Q)$
is infinite. Hence a union of finitely many $B$-orbits can not
be a $G$-orbit. It follows that $G$ acts transitively on $\Omega$.
\end{proof}
\begin{corollary}
Let $k$ be a cubic number field, and $\Gamma_1,\Gamma_2$ lattices
in $\C$ which are contained in $k$.
Then the two elliptic curves $\C/\Gamma_1$ and $\C/\Gamma_2$
are isogenous.
\end{corollary}
\begin{proof}
This follows because two lattice $\Gamma_1$, $\Gamma_2$ in $\C$ have
isogenous quotient elliptic curves if and only if there is
an element $A\in PGL_2(\Q)$ such that the associated fractional
linear transformation carries $\Gamma_1$ to $\Gamma_2$.
\end{proof}

\begin{corollary}
Let $k$ be a cubic number field,
$\tau$ a primitive element for the field extension
$k/\Q$ and $\Gamma$ a lattice in $\C^2$ which is
contained in $k^2$.
Let $E$ be an elliptic curve for which there exists a surjective
holomorphic Lie group homomorphism $\pi:\C^2/\Gamma\to E$.

Then $E$ is isogenous to $\C/\left<1,\tau\right>_\Z$.
\end{corollary}
\begin{proof}
Such a surjective homomorphism $\pi$ is induced by a linear
map $\tilde\pi:\C^2\to\C$. Since $\Gamma$ has $\Z$-rank $3$, there
is a non-zero element of $\Gamma$ in the kernel of $\tilde\pi$.
Therefore there are $\mu\in\C^*$ and $a_i\in k$ such that
$\tilde\pi(x_1,x_2)=\mu(a_1x_1+a_2x_2)$. It follows that
\[
\frac{1}{\mu}\tilde\pi(\Gamma)\subset k.
\]
Thus $E$ is isogenous to $\C/\left<1,\tau\right>_\Z$
by the preceding claim.
\end{proof}

\begin{proof}[Proof of the theorem]
For $R\in\{\Q,\R,\C\}$ the $R$-module generated by
$\Gamma$ will be denoted by $\Gamma_R$. Let $H=\Gamma_\R\cap i\Gamma_\R$.
Then for every $\gamma\in\Gamma\setminus H$ we obtain a surjective
Lie group homomorphism onto an elliptic curve as follows:
Let $L_\gamma$ denote the quotient map from $\C^2$ to $Q=\C^2/\C\cdot\gamma$.
By construction the image $L_\gamma(\Gamma)$ has rank at most two
(because $\Z\cdot\gamma\subset\ker L_\gamma$).
On the other hand $\gamma\not\in H$ implies $L_\gamma(H)=Q$ and
therefore $L_\gamma(\Gamma_\R)=Q$. It follows that $L_\gamma(\Gamma)$
is a lattice in $Q$ and $Q/L_\gamma(\Gamma)$ is an elliptic curve.

By choosing a basis for the complex vectorspace $\C^2$ inside $\Gamma\setminus H$ we may assume that
\[
\Gamma = \left< \begin{pmatrix} 1 \\ 0 \end{pmatrix}
\begin{pmatrix} 0 \\ 1 \end{pmatrix}
\begin{pmatrix} \alpha \\ \beta \end{pmatrix}
\right>_\Z
\]
with $\alpha,\beta\in\C\setminus\R$.

Consider
\[
\gamma=\gamma_{m,n,p}=n \begin{pmatrix} 1 \\ 0 \end{pmatrix}
+m\begin{pmatrix} 0 \\ 1 \end{pmatrix}
+p\begin{pmatrix} \alpha \\ \beta \end{pmatrix}
\]
for $(m,n,p)\in\Q^3$.
Let $I$ denote the $\Q$ vector subspace of $\Q^3$ of all
$(m,n,p)$ for which $\gamma_{m,n,p}\in H$.
Evidently $\dim_\Q(I)\in\{0,1,2\}$.

Let us first deal with the case $\dim_\Q(I)=2$.
Then $H\cap\Gamma$ has $\Z$-rank two and $H/(H\cap\Gamma)$ is compact.
Consider the linear projection $\pi:\C^2\to\C^2/H$.
Then $\pi(\Gamma)$ is a subgroup of $(\C,+)$ of $\Z$-rank one.
Hence $\pi(\Gamma)\simeq\Z$. It follows that
$T\simeq\C^*\times H/(H\cap\Gamma)$.

Thus we may assume that $\dim_\Q(I)\in\{0,1\}$.

We claim that this implies 
$\left<1,\alpha\right>_\Q\ne\left<1,\beta\right>_\Q$.
Indeed, if $\alpha=x\beta+y$ for some $x,y\in\Q$,
then 
\[
\begin{pmatrix} x\beta \\ \beta \end{pmatrix}
=
\begin{pmatrix} x\beta+y\\ \beta \end{pmatrix}
-
\begin{pmatrix} y \\ 0 \end{pmatrix}
\in \Gamma_\Q.
\]
Consequently $H=\{(xt,t):t\in\C\}$ and $\rank_\Z(H\cap\Gamma)=\dim_\Q(I)=2$.

For $(m,n,p)\in\Q^3\setminus\{(0,0,0)\}$ we define $E_{m,n,p}$ as the 
$\Q$-vector space
\[
E_{m,n,p}=\left< 
m+p\beta,n+p\alpha,m\alpha-n\beta
\right>_\Q.
\]
One verifies easily that for all  $(m,n,p)\in\Q^3\setminus\{(0,0,0)\}$ the matrix
\[
\begin{pmatrix} m & 0 & p \\ n & p & 0 \\ 0 & m & -n \end{pmatrix}
\]
has rank $2$ which implies that $E_{m,n,p}$ is a $\Q$-plane in $V=\left<1,\alpha,\beta\right>_\Q$.

Thus $(m,n,p)\mapsto E_{m,n,p}$ defines a map from $\Q^3\setminus\{(0,0,0)\}$
to the Grassmann variety $G_\Q$ of $\Q$-planes in $V$.

If $(m,n,p)\in\Q^3\setminus I$ we have a projection 
$L_{m,n,p}:\C^2\to\C^2/\C\gamma_{m,n,p}$
which in coordinates can be described as
\[
(z_1,z_2)\mapsto (m+p\beta)z_1 - (n+p\alpha)z_2.
\]
In this case $E_{m,n,p}$ can be identified with the image of $\Gamma_\Q$
under the projection $L_{m,n,p}$.

Let $E$ be an arbitrary $\Q$-plane in $V$. 
There are two possibilities: Either $1\in E$ or $1\not\in E$.
In the first case
$E=\left<1,m\alpha-n\beta\right>_\Q$ for some $(m,n)\in\Z^2\setminus\{(0,0)\}$.
Then $E=\psi(m,n,0)$.
Let us now discuss the second case.
Since $E$ is a $\Q$-hyperplane in $V$, we have that both
$E\cap\left<1,\alpha\right>_\Q$ and $E\cap\left<1,\beta\right>_\Q$
have $\Q$-dimension one.
Therefore $E$ is the direct sum of two $\Q$-lines which arise as intersections
of $E$ with $\left<1,\alpha\right>_\Q$ resp.~$\left<1,\beta\right>_\Q$.
Thus
\[
E=\left<m+p\alpha,n+p\beta\right>_\Q
\]
for some $(m,n,p)\in\Q^3\setminus\{(0,0,0)\}$.
It follows that $E=\psi(m,n,p)$.

In this way we have shown that the map $\psi:
\Q^3\setminus\{0\}\to G_\Q$ is surjective.

Observe that $I$ is at most one-dimensional
and that $\psi(m,n,p)$ depends only on $[m:n:p]$.
Therefore we see: {\em With at most one exception every element of
$G_\Q$ is in the image of $\Q^3\setminus I$ under $\psi$.}
Now let $\Zeta$ denote the set of all $k\in\Z$ for which
$\left<1,\beta+k\alpha\right>_\Q$ is contained in
the image $\psi(\Q^3\setminus I)$.
Then $\Z\setminus\Zeta$ contains at most one element.

Now  $\C/\left<1,\alpha\right>_\Z$
is isogenous to  $\C/\left<1,\beta+k\alpha\right>_\Z$
for all $k\in\Zeta$.

Thus for every $k\in\Zeta$ there is a 
matrix
\[
A_k=\begin{pmatrix} a_k & b_k \\ c_k & d_k 
    \end{pmatrix}
\in GL_2(\Q)
\]
such that $A_k(\alpha)=\frac{a_k\alpha+b_k}{c_k\alpha+d_k}=\beta+k\alpha$.
Then $A_k(\alpha)-A_n(\alpha)=(k-n)\alpha$ for all $k,n\in\Zeta$.
It follows that
\[
(a_k\alpha+b_k)(c_n\alpha+d_n) -
(a_n\alpha+b_n)(c_k\alpha+d_k) 
=(k-n)\alpha(c_k\alpha+d_k) (c_n\alpha+d_n)
\]
Thus for all $(k,n)\in\Zeta\times\Zeta$ there is a $\Q$-polynomial $P_{k,n}$
of degree at most two
such that
\[
\alpha^3(k-n)c_kc_n=P_{k,n}(\alpha).
\]
Therefore either $\deg(\alpha)\le 3$
or $(k-n)c_kc_n=0$ for all $k,n\in\Zeta$.

However, $A_k(\alpha)=\beta+k\alpha$ combined with $\left<1,\alpha\right>_\Q
\ne\left<1,\beta\right>_\Q$ implies that $c_k\ne 0$.

Therefore $\deg(\alpha)\le 3$.

Furthermore $A_k(\alpha)=\beta+k\alpha$ implies $\beta\in\Q(\alpha)$.

For this reason we may deduce that
 $\Gamma\subset k^2$ for some number field $k$ of degree
at most three.

Finally we recall that $T$ splits into a direct product of $\C^*$ and an
elliptic curve if $\deg(k/\Q)=2$  (lemma~\ref{q-field}).
\end{proof}

\section{Isogeny vs. isomorphism}

\begin{lemma}
Let $z,w$ be complex numbers with $\Im(z),\Im(w)>1$.

Then the elliptic curves $\C/\left<1,z\right>_{\Z}$ and
$\C/\left<1,w\right>_{\Z}$ are biholomorphic if and only if
$z-w\in\Z$.
\end{lemma}
\begin{proof}
This follows easily from the well-known fact that
\[
F=\left\{z\in\C: |z|>1, \Im(z)>0 \text{ and }|\Re(z)|<\frac{1}{2} \right\}
\]
is a fundamental domain for the $PSL_2(\Z)$-action on the upper halfplane.
\end{proof}

\begin{proposition}
Let $\Gamma$ be a discrete subgroup of $\C^2$ of $\Z$-rank three.

Then there exists non-isomorphic quotient elliptic curves.
\end{proposition}
\begin{proof}
Without loss of generality we may assume that
\[
\Gamma=\left< \begin{pmatrix} 1 \\ 0 \end{pmatrix},
              \begin{pmatrix} 0 \\ 1 \end{pmatrix},
              \begin{pmatrix} \alpha \\ \tau \end{pmatrix}
\right>_{\Z}
\]
with $\Im(\tau)>0$ and $\alpha\in\C$.
For every $m\in\Z$ the group $\Gamma$ contains $(m,1)$. The quotient
by the complex line through this element is given by
\[
(z_1,z_2)\mapsto mz_2-z_1
\]
The image of $\Gamma$ in $\C$ under this projection is
\[
\Lambda_m=\left< -1, m, m\tau-\alpha\right>_\Z= 
\left<1, m\tau-\alpha\right>_\Z.
\]
By the preceding lemma the quotients $\C/\Lambda_m$ and $\C/\Lambda_n$
are not isomorphic for integers $m,n$ with
\[
m > n > \frac{1+\Im(\alpha)}{\Im(\tau)}
\]
because $n > \frac{1+\Im(\alpha)}{\Im(\tau)}$ is equivalent to
$\Im(n\tau-\alpha)>1$ and $m\tau-n\tau$ is never contained in $\Z$.
\end{proof}
\begin{example}
Let $\Gamma=\Z\times\Z[i]\subset\C\times\C$.
Then $(z_1,z_2)\mapsto z_1-mz_2$ maps $\C^2/\Gamma$ onto the elliptic
curve $\C/\left<1,mi\right>_{\Z}$.
\end{example}


\begin{thebibliography}{Bla}
\bibitem{C}
Cousin, P.:
Sur les fonctions triplement p\'eriodiques de deux variables.
\sl Acta Math. \bf 33, \rm 105--232 (1905)

\end{thebibliography}
\end{document}